\documentclass[12pt,reqno]{amsart}
\usepackage{amsmath, amsthm, amssymb}
\usepackage{mathtools}
\usepackage{amsaddr}
\usepackage[english]{babel}

\usepackage[margin=1in]{geometry}
\usepackage{parskip}
\usepackage[shortlabels]{enumitem}

\usepackage{xcolor}

\usepackage[bookmarks,colorlinks,breaklinks]{hyperref}
\hypersetup{linkcolor=red,citecolor=blue,
  filecolor=dullmagenta,urlcolor=blue}

\usepackage[numbers]{natbib}
\usepackage{cleveref}

\usepackage{graphicx}
\usepackage{tikz}
\usetikzlibrary{matrix, arrows}
\usetikzlibrary{cd}

\makeatletter
\def\thm@space@setup{%
  \thm@preskip=0.5em\thm@postskip=\thm@preskip%
}
\makeatother

\newtheoremstyle{named}{}{}{\\itshape}{}{\bfseries}{.}{.5em}{\thmnote{#3's }#1}
\theoremstyle{named}

\theoremstyle{plain}
\newtheorem{thm}{Theorem}[section]

\newtheorem{cor}[thm]{Corollary}
\theoremstyle{definition}

\theoremstyle{remark}
\newtheorem{rmk}[thm]{Remark}
\crefformat{footnote}{#2\footnotemark[#1]#3}

\newcommand{\mathfont}{\mathbf}

\newcommand{\ZZ}{\mathfont Z}
\newcommand{\Zhat}{\widehat{\ZZ}}

\newcommand{\QQ}{\mathfont Q}
\newcommand{\Qbar}{\overline{\QQ}}

\newcommand{\F}{\mathfont F}
\newcommand{\FF}{\mathfont F}
\newcommand{\Fbar}{\overline{\FF}}

\newcommand{\fG}{\mathfrak{G}}

\newcommand{\fg}{\mathfrak{g}}

\newcommand{\bA}{\mathfont{A}}

\newcommand{\bP}{\mathfont{P}}

\newcommand{\cO}{\mathcal{O}}

\DeclareFontFamily{OT1}{rsfs}{}
\DeclareFontShape{OT1}{rsfs}{n}{it}{<-> rsfs10}{}
\DeclareMathAlphabet{\mathscr}{OT1}{rsfs}{n}{it}

\newcommand{\onto}{\twoheadrightarrow}

\DeclareMathOperator{\Hom}{Hom}
\DeclareMathOperator{\Ker}{ker}

\DeclareMathOperator{\Gal}{Gal}

\DeclareMathOperator{\der}{der}

\DeclareMathOperator{\Frob}{Frob}

\newcommand{\Zl}{\mathfont{Z}_{\ell}}

\newcommand{\rhobar}{{\overline{\rho}}}

\DeclareMathOperator{\GL}{GL}
\DeclareMathOperator{\SL}{SL}

\DeclareMathOperator{\Image}{Im}

\DeclareMathOperator{\Lift}{Lift}

\DeclareMathOperator{\mr}{m.r.}

\makeatother
\begin{document}


\title{Remarks on the inverse Galois problem over function fields \\
}

\author{Shiang Tang}


\email{shiangtang1989@gmail.com}


\maketitle

\begin{abstract}
In this paper, we prove new instances of the inverse Galois problem over global function fields for finite groups of Lie type. This is done by constructing compatible systems of $\ell$-adic Galois representations valued in a semisimple group $G$ using Galois theoretic and automorphic methods, and then proving that the Galois images are maximal for a set of primes of positive density using a classical result of Larsen on Galois images for compatible sytems.  
\end{abstract}

\section{Introduction}
The classical inverse Galois problem asks if every finite group is isomorphic to the Galois group of a finite Galois extension of $\QQ$. Even though it is still far from being completely solved, there are many established and beautiful approaches to it, see \cite{mm:igp}. In the past two decades there has been new progress on the inverse Galois problem via the theories of motives and automorphic forms, for example, \cite{kls}, \cite{yun}, \cite{boxer:e6}, and \cite{tan:spin}. To give a simple case of interest, let us consider an elliptic curve $E$ over $\QQ$, its $\ell$-adic Tate modules for varying primes $\ell$ give rise to a compatible system of $\ell$-adic Galois representations $\rho_{\ell} \colon \Gal(\Qbar/\QQ) \to \GL_2(\Zl)$, where `compatible' means that the characteristic polynomial of the image of the Frobenius element at any unramified prime has $\QQ$-coefficients, and it is independent of $\ell$. By Serre's open image theorem, $\rho_{\ell}$ is surjective for all but finitely many $\ell$; in particular, $\GL_2(\FF_{\ell})$ is isomorphic to the Galois group of a finite Galois extension of $\QQ$. 

In this paper, we consider the inverse Galois problem over a global function field and prove new instances of it by constructing compatible systems of Galois representations valued in a semisimple algebraic group. Our main theorem is the following:

\begin{thm} \label{thm:intromain}
Let $F$ be a global function field of characteristic $p$ (e.g. $F=\F_p(T)$). 
Let $G$ be a split semisimple group. Then there is a set $\Sigma$ of rational primes of positive Dirichlet density, such that for $\ell \in \Sigma$, $G(\F_{\ell})$ is isomorphic to the Galois group of a finite Galois extension of $F$. 
\end{thm}

We remark that the inverse Galois problem is still open for global function fields, and indeed for every global field. 
We also remark that the inverse Galois problem is known for $\Fbar_p(T)$ (which is not a global field) due to a well-known result of Harbater and Raynaud.
Our theorem is in contrast to the following theorem due to Fried--V\"olklein, Jarden and Pop \cite[Proposition 3.3.9]{harbater-patching}, proven using a completely different approach:

\begin{thm}
Let $\fG$ be a finite group. Then for all but finitely many finite fields $\FF$, $\fG$ is isomorphic to the Galois group of a finite Galois extension of $\FF(T)$. 
\end{thm}

To the best of our knowledge, the above theorems are the only known results on the inverse Galois problem over global function fields.

This paper is organized as follows. In Section 2, we prove that if $F$ is a global function field and $G$ is a semisimple algebraic group, then any continuous representation of $\Gal(\overline{F}/F)$ into $G$ in characteristic $\ell$ satisfying certain technical conditions lifts to an $\ell$-adic representation that is part of a compatible system of Galois representations. This follows from a lifting theorem of Fakhruddin--Khare--Patrikis and work of L. Lafforgue on the Langlands correspondence for function fields. In Section 3, we construct Galois representations in characteristic $\ell$ to which the results in Section 2 apply. Then we show that the resulting compatible system of $\ell$-adic Galois representations has maximal Galois images for a positive proportion of primes $\ell$ using a result of Larsen on maximality of Galois images for compatible systems. 

\emph{Acknowledgements}: I would like to thank Jeremy Booher and Anwesh Ray for helpful conversations and comments. I also thank the anonymous referee for helpful comments and corrections. 


\subsection{Notation}
Let $G$ be a split connected reductive group defined over $\Zl$. We let $G^{\der}$ denote the derived subgroup of $G$, and we denote by $\fg^{\der}$ the Lie algebra of $G^{\der}$. 
Let $G^{ab}$ be the quotient $G/G^{\der}$.
Let $k$ be a finite extension of $\FF_{\ell}$ and let $\cO=W(k)$ be its ring of Witt vectors. Let $F$ be a global function field of characteristic $p$. Fix a separable closure $\overline{F}$ of $F$. Let $\Gamma_F=\Gal(\overline{F}/F)$ be the absolute Galois group of $F$, and for a finite set $S$ of places of $F$, let $\Gamma_{F,S}$ be the Galois group of the maximal extension of $F$ inside $\overline{F}$ that is unramified outside the primes in $S$. Given a continuous representation $\rhobar \colon \Gamma_F \to G(k)$, let $\bar{\mu} \colon \Gamma_F \to G^{ab}(k)$ be its pushforward via the natural map $\nu \colon G \to G^{ab}=G/G^{\der}$. We write $\rhobar(\fg^{\der})$ for $\fg^{\der}_k$ equipped with an action of $\Gamma_F$ via the composite of $\rhobar$ and the adjoint representation. We write $\rhobar(\fg^{\der})^*$ for the Tate dual $\Hom(\rhobar(\fg^{\der}),\mu_{\ell})$, which is isomorphic to $\rhobar(\fg^{\der})(1)$, the twist of $\rhobar(\fg^{\der})$ by the mod-$\ell$ cyclotomic character. 

\section{Lifting Galois representations over function fields} \label{sec:lifting} 
We begin with the following variant of a lifting theorem of Fakhruddin--Khare--Patrikis \cite[Theorem E]{fkp:aut} in the function field case. Recall that $F$ is a function field of characteristic $p$, let $\ell$ be a prime distinct from $p$. 

\begin{thm} \label{thm:fkp}
Suppose $\ell \gg_G 0$ (i.e. $\ell$ is large enough relative to the root datum of $G$). Let 
$\rhobar \colon \Gamma_{F,S} \to G(k)$ be a continuous representation. Fix an integer $t$, and fix a lift $\mu \colon \Gamma_F \to G^{ab}(\cO)$ of $\bar{\mu} \colon \Gamma_F \to G^{ab}(k)$. Assume that
\begin{enumerate}
    \item $H^1(G(K/F), \rhobar(\fg^{\der})^*)=0$, where $K=F(\rhobar, \mu_{\ell})$.
    \item $H^0(\Gamma_F, \rhobar(\fg^{\der}))=H^0(\Gamma_F, \rhobar(\fg^{\der})^*)=0$.
    \item There is no surjection of $\F_{\ell}[\Gamma_F]$-modules $\rhobar(\fg^{\der}) \onto W$ for some $\F_{\ell}[\Gamma_F]$-module subquotient $W$ of $\rhobar(\fg^{\der})^*$.
\end{enumerate}
Then there exist a finite set of primes $\widetilde{S}$ containing $S$ and a lift $\rho \colon \Gamma_{F,\widetilde{S}} \to G(\cO)$ such that
\begin{itemize}
    \item $\rho$ has multiplier $\mu$.
    \item $\rho(\Gamma_F) \supset \widehat{G^{\der}}(\cO) \coloneqq \Ker( G^{\der}(\cO) \to G^{\der}(k) )$.
    \item For any $v \in S$, $\rho|_{\Gamma_{F_v}}$ is an $\cO$-point of $R_{\rhobar|_{\Gamma_{F_v}}}^{\mr,\Box}$, where $R_{\rhobar|_{\Gamma_{F_v}}}^{\mr,\Box}$ is an $\cO$-formally smooth quotient of the universal lifting ring for $\rhobar|_{\Gamma_{F_v}}$ constructed in \cite{bct}. 
\end{itemize}
\end{thm}

\begin{proof}
We apply \cite[Theorem E]{fkp:aut}. Its assumptions are clearly met except the existence of local lifts. But this follows immediately from \cite[Theorem 1.1]{bct}, which shows that the universal lifting ring of $\rhobar|_{\Gamma_{F_v}}$ for $v \in S$ has an $\cO$-formally smooth quotient $R_{\rhobar|_{\Gamma_{F_v}}}^{\mr,\Box}$. We may and do impose the corresponding formally smooth local deformation condition $\Lift_{\rhobar|_{\Gamma_{F_v}}}^{\mr}$ in the global deformation argument in \cite{fkp:aut}, and the theorem follows. Note that by \cite[Remark 1.3]{fkp} one does not need to enlarge the coefficient ring $\cO$ since a formally smooth deformation condition exists at every finite place.
\end{proof}

\begin{rmk}
Note that we have removed the assumption on the existence of local lifts in \cite[Theorem E]{fkp:aut}, and we do not require extending the coefficient ring $\cO$ as opposed to \emph{loc. cit.}
\end{rmk}

\begin{rmk}
For $\ell \gg_G 0$, the assumptions in Theorem \ref{thm:fkp} are satisfied if 
\begin{itemize}
    \item $\rhobar|_{F(\mu_{\ell})}$ is absolutely irreducible, and
    \item $[F(\mu_{\ell}):F]>a_G$ for a constant $a_G$ depending only on the root datum of $G$.
\end{itemize}
In fact, this follows from the proof of \cite[Corollary A.7]{fkp}.
\end{rmk}

\begin{rmk}
The method for proving Theorem \ref{thm:fkp} is purely Galois-theoretical, which uses only Galois cohomology and arithmetic duality theorems. One needs to allow additional ramifications in order to annihilate the global obstructions contained in the Selmer groups. In contrast, using automorphic techniques, one can obtain similar lifting theorem \emph{without} adding ramifications, assuming that the residual representation is irreducible, see \cite[Theorem 5.6]{bfhkt} (resp. \cite[Theorem 4.3.1]{blggt}) for a theorem of this kind for $G=\GL_n$ in the case of function fields (resp. number fields). However, even in the $\GL_n$ case, Theorem \ref{thm:fkp} applies to a much larger class of residual representations $\rhobar$, including reducible representations, at the expense of adding ramifications. 
One can prove similar theorems for general $G$ without enlarging the ramification locus by comparing the $\GL$-deformation ring and $G$-deformation ring, assuming that the restriction of the residual representation is irreducible after being composed with a faithful representation $G \to \GL(V)$. We do not pursue it here, see \cite[Theorem 3.4]{pt:gspin} a theorem of this kind for the odd spin groups.
\end{rmk}

We now deduce the following result, which will be used in proving Theorem \ref{thm:intromain}. For generalities on compatible systems of Galois representations valued in a general reductive group, see \cite[Section 6]{bhkt}. 

\begin{cor} \label{cor:comp system}
Retain the assumptions of Theorem \ref{thm:fkp}. Suppose moreover that $G$ is semisimple. Then there exist a finite set of places $\widetilde{S}$ containing $S$, and
a compatible system of $G$-valued Galois representations (in the sense of \cite[Definition 6.1]{bhkt}) 
$\rho_{\lambda} \colon \Gamma_{F,\widetilde{S}} \to G(\Qbar_{\lambda})$ indexed by the prime-to-$p$ places $\lambda$ of $\Qbar$ such that for some $\lambda|\ell$, $\rho_{\lambda}$ lifts $\rhobar$ and the image of $\rho_{\lambda}$ is Zariski-dense for all places $\lambda$ of $\Qbar$. 
\end{cor}

\begin{proof}
By Theorem \ref{thm:fkp}, for a finite set of places $\widetilde{S}$ containing $S$, $\rhobar$ lifts to a continuous representation $\rho \colon \Gamma_{F,\widetilde{S}} \to G(\Qbar_{\ell})$ 
such that $\Image(\rho) \supset \widehat{G^{\der}}(\cO)$; in particular, since $G$ is semisimple, this implies that $\rho$ has Zariski-dense image. 
By \cite[Theorem 6.5]{bhkt}, $\rho$ belongs to a compatible system of $G$-valued Galois representations claimed above. 
\end{proof}

\begin{rmk}
Note that in the proof of \cite[Theorem 6.5]{bhkt}, the semisimplicity assumption on $G$ is needed to apply Lafforgue's theorem \cite[Theoreme VII.6]{lafl} after composing $\rho$ with a faithful representation of $G$. 
\end{rmk}

\section{The inverse Galois problem for function fields} \label{sec:igp}

In this section we prove Theorem \ref{thm:intromain}. 


\begin{proof}[Proof of Theorem \ref{thm:intromain}]
We may and do assume that $G$ is simply-connected. 
We will construct a suitable residual representation $\rhobar \colon \Gamma_F \to G(k)$, deform it to an $\ell$-adic representation belonging to a compatible system using Corollary \ref{cor:comp system}, then apply a result of Larsen to deduce that the compatible system has full image $G(\Zl)$ for a positive density set of primes $\ell$. 

Arguing as in the first two paragraphs of the proof of \cite[Lemma 10.12]{bhkt}, we may construct a continuous representation
$\rhobar \colon \Gamma_{F} \to G(\overline{\F_{\ell}})$ such that
\begin{itemize}
    \item The image of $\rhobar$ has prime-to-$\ell$ order.
    \item $\rhobar$ is absolutely irreducible.
    \item The fields $F(\rhobar)$ and $F(\mu_{\ell})$ are linearly disjoint over $F$. 
\end{itemize}
More precisely, we take $\rhobar$ to be the mod-$\ell$ reduction of the Coxeter homomorphism constructed in \emph{loc.cit} ($\ell$ to be specified later).  
Since this Coxeter homomorphism is depends only on the field $F$ and the group $G$, the first point is met as long as $\ell$ is large enough. 
The second point then follows from \cite[Proposition 10.10, (ii)]{bhkt}.  
For the third point, it suffices to choose $\ell$ such that the order of $p$ modulo $\ell$ is prime to $d := [F(\rhobar):F]$ (which is independent of $\ell$). It is easy to see that there are infinitely many such primes $\ell$: for example, consider the prime factors of $p^N - 1$ where $(N,d)=1$. For any such $\ell$, it follows that $[F(\mu_{\ell}) : F]$ is prime to $d$, whence the linear disjointness. 

We now check the assumptions (1)-(3) in Theorem \ref{thm:fkp}. Recall that $\ell \gg_G 0$. We also assume that $\ell \nmid p-1$. 
By the first point, $\Gal(K/F)$ has prime-to-$\ell$ order, which implies (1). (2) follows from absolute irreducibility and \cite[Lemma A.2]{fkp}. For (3), let $\{W_i\}_{i \in I}$ be the irreducible summands of $\rhobar(\fg)$ (note that $\rhobar(\fg)$ is semisimple). If $\rhobar(\fg)$ and $\rhobar(\fg)^* \cong \rhobar(\fg)(1)$ had a common constituent, there would be an isomorphism $W_i \cong W_j(1)$ for some $i,j \in I$. By the third point, we can choose $\sigma \in \Gamma_F$ acting trivially on $W_i$ and $W_j$ but non-trivially on $F(\mu_{\ell})$, which contradicts to that $W_i \cong W_j(1)$. 
By Corollary \ref{cor:comp system}, there is a compatible system $\rho_{\lambda} \colon \Gamma_{F} \to G(\Qbar_{\lambda})$ indexed by the prime-to-$p$ places $\lambda$ of $\Qbar$ with Zariski-dense images. By \cite[Proposition 6.6]{bhkt} (which is an application of \cite[Theorem 3.17]{larsen}), after passing to an equivalent compatible system, there exists a set $\Sigma$ of rational primes of positive Dirichlet density, such that if $\ell \in \Sigma$ and $\lambda$ is a place of $\Qbar$ above $\ell$, then $\rho_{\lambda}$ has image equal to $G(\Zl)$. In particular, $G(\F_{\ell})$ is isomorphic to the Galois group of a finite Galois extension of $F$. 
\end{proof}

\begin{rmk}
In the case when $F=\F_p(T)$, Theorem \ref{thm:intromain} answers instances of Abhyankar's total arithmetical conjecture \cite[Conjecture 5.3]{hops}. On the other hand, for a smooth projective curve $X$ defined over $\Fbar_p$, Ahbyankar's conjecture (now a theorem of Harbater and Raynaud, see \cite{harbater} and \cite{raynaud}) implies that every finite group arises as the Galois group of a branched cover of $X$. 
\end{rmk}


\begin{rmk}
If one considers restricting the ramification locus, then the conclusion of Theorem \ref{thm:intromain} may be false. As a trivial example, $\pi_1(\bP^1_{\FF_p}) \cong \Zhat$. For a nontrivial example, consider $\bA^1_{\FF_p}$. By \cite[Conjecture 5.2]{hops}, a finite group $A$ is a quotient of $\pi_1(\bA^1_{\FF_p})$ if and only if $A$ is an extension of a cyclic group by a quasi-$p$ group. In particular, a finite simple group is a quotient of $\pi_1(\bA^1_{\FF_p})$ if and only if its order is divisible by $p$. 
\end{rmk}

\begin{rmk}
The proof of Theorem \ref{thm:intromain} in fact shows that the statement holds for $G(\ZZ/\ell^n \ZZ)$ as well. 
\end{rmk}


\end{document}